\documentclass[11pt]{article}
\usepackage{amsmath}
\usepackage{amsthm}
\usepackage{amsfonts}
\usepackage{amssymb}
\usepackage{latexsym}

\usepackage[colorlinks,
linkcolor=blue,
anchorcolor=blue,
citecolor=blue,
urlcolor=blue
]{hyperref}

\hypersetup{
	pdftitle={Descent generating polynomials for (n-3)- and (n-4)-stack-sortable (pattern-avoiding) permutations},
	pdfauthor={Philip B. Zhang and Sergey Kitaev},
	pdfsubject={Published in Discrete Applied Mathematics 372 (2025), 1-14},
	pdfkeywords={stack-sortable permutations; descent statistic; Eulerian polynomial; pattern avoidance}
}

\newtheorem{thm}{Theorem}
\newtheorem{lem}[thm]{Lemma}

\newtheorem{conj}[thm]{Conjecture}
\newtheorem{rem}[thm]{Remark}

\newcommand{\sortable}[2]{\ensuremath{\mathcal{W}_{#1,#2}}}
\newcommand{\exactly}[2]{\ensuremath{\mathcal{E}_{#1,#2}}}

\newcommand{\sym}{\mathfrak{S}}

\newcommand{\st}{\mbox{\ensuremath{\ast}}}
\newcommand{\q}{\mbox{\ensuremath{\scriptstyle{?}}}}

\newcommand{\glob}[1]{\ensuremath{\langle\, #1\,\rangle}}

\DeclareMathOperator{\des}{des}

\begin{document}

\begin{center}
	{\large \bf  Descent generating polynomials for ($n-3$)- and ($n-4$)-stack-sortable (pattern-avoiding) permutations}
\end{center}

\begin{center}
	Philip B. Zhang$^{a}$ and Sergey Kitaev$^{b}$
	\\[6pt]
	
	$^{a}$School of Mathematical Sciences, Tianjin Normal University,\\
	Tianjin 300387, P. R. China\\[6pt]
	
	$^{b}$Department of Mathematics and Statistics, University of Strathclyde,\\
	26 Richmond Street, Glasgow G1 1XH, United Kingdom\\[6pt]
	
	Email: $^{a}${\tt zhang@tjnu.edu.cn},
	$^{b}${\tt sergey.kitaev@strath.ac.uk}
\end{center}

\noindent\textbf{Abstract.}
In this paper, we find distribution of descents over $(n-3)$- and $(n-4)$-stack-sortable permutations in terms of Eulerian polynomials. Our results generalize the enumeration results by Claesson, Dukes, and Steingr\'{\i}msson  on $(n-3)$- and $(n-4)$-stack-sortable permutations. Moreover, we find distribution of descents on $(n-2)$-, $(n-3)$- and $(n-4)$-stack-sortable permutations that avoid any given pattern of length 3, which extends known results in the literature on distribution of descents over pattern-avoiding 1- and 2-stack-sortable permutations. Our distribution results also give enumeration of $(n-2)$-, $(n-3)$- and $(n-4)$-stack-sortable permutations avoiding any pattern of length 3. One of our conjectures links our work to stack-sorting with restricted stacks, and the other conjecture states that 213-avoiding permutations sortable with $t$ stacks are equinumerous with 321-avoiding permutations sortable with $t$ stacks for any $t$.

\noindent {\bf AMS Classification 2010:} 05A15


\section{Introduction}

%The main purpose of this paper is to give  two formulas of  the generating functions of the descent statistic for  ($n-3$)- and ($n-4$)-stack-sortable permutations  in terms of Eulerian polynomials.

Let $n$ be  a positive integer. A permutation of length $n$ is a rearrangement of the set $[n]:=\{1,2,\ldots,n\}$ and $\varepsilon$ denotes the empty permutation.
For a permutation $\pi=\pi_1\pi_2\cdots \pi_n$ with $\pi_i=n$,  the {\em stack-sorting operation} $\mathcal{S}$  is recursively defined as follows, where $\mathcal{S}(\varepsilon) =\varepsilon$:
$$\mathcal{S}(\pi) = \mathcal{S}(\pi_1\cdots \pi_{i-1})\,  \mathcal{S}(\pi_{i+1}\cdots \pi_n)\, n.$$
%It is well known that
A permutation $\pi$ is {\em $t$-stack-sortable} if the permutation $\mathcal{S}^t(\pi)$, obtained by applying the stack-sorting operation $t$ times,  is the identity permutation. The study of ($t$-)stack-sortable permutations has  attracted great attention in mathematics and theoretical computer science (see, e.g.\ \cite{BBS2010,Bona2002Symmetry,BKLPRW,Cerbai2021,Cerbai2020Stack,Claesson2010Permutations,Defant2021Asymptotics,Defant2021Stack,Zhang2015real}).

Clearly, all permutations are $(n-1)$-stack-sortable. Knuth~\cite{Knuth1969Art} discovered that  the Catalan numbers $C_n=\frac{1}{n+1}\binom{2n}{n}$ count the number of permutations of length $n$ whose resulting permutations  are 1-stack-sortable permutations of length $n$. It turns out that 1-stack-sortable permutations are precisely  231-avoiding permutations, where a permutation $\pi_1\pi_2\cdots\pi_n$ avoids a pattern $p=p_1p_2\cdots p_k$ (which is also a permutation) if there is no subsequence $\pi_{i_1}\pi_{i_2}\cdots\pi_{i_k}$ such that $\pi_{i_j}<\pi_{i_m}$ if and only if $p_j<p_m$ \cite{Kitaev2011Patterns}. For the cases of $t=2$ and $t=n-2$, enumeration results for $t$-stack-sortable permutations were given by West~\cite{West1990Permutations},
and  for $t=n-3$ and $t=n-4$ by
Claesson,  Dukes, and Steingr\'\i msson \cite{Claesson2010Permutations}.
For an overview of properties of $t$-stack-sortable permutations,
see \cite{Bona2002Symmetry,Knuth1969Art,West1990Permutations} and \cite[Sec 2.1]{Kitaev2011Patterns}. More recent results related to stack-sortings can be found here \cite{BKLPRW,Cerbai2020Stack,Defant2022Troupes, Defant2021Asymptotics,  Defant2021Stack}.

In this paper, we consider the descent generating polynomials for (pattern-avoiding) $t$-stack-sortable permutations for certain values of $t$.
Given a permutation $\pi=\pi_1 \pi_2 \cdots \pi_n$,  the {\em descent statistic } on $\pi$ is defined as
the number of $i\in [n-1]$ such that $\pi_i > \pi_{i+1}$.
Denote by  $\mathfrak{S}_n$  the set of permutations of $[n]$.
The {\em Eulerian polynomial} $A_n(x)$ is defined as the  generating function  of the descent statistic over $\mathfrak{S}_n$, namely,
\begin{align*}
	A_n(x) := \sum_{\pi \in \mathfrak{S}_n}x^{\des (\pi)}.
\end{align*}
It is known that the Eulerian polynomials can be computed inductively by
$$A_0(x)=1, \ \ A_n(x)=\sum_{k=0}^{n-1}{n\choose k}A_k(x)(x-1)^{n-1-k},\ \ n\geq 1.$$
Also, for $n\geq 1$, the Eulerian polynomials can be computed as
$$A_n(x)=\sum_{k=1}^{n}k!S(n,k)(x-1)^{n-k}$$
where $S(n,k)$ is the {\em Stirling numbers of the second kind}.

Denote by $\sortable{n}{t}$ the set of $t$-stack-sortable permutations of $[n]$. Let
\begin{align*}
	W_{n,t}(x) := \sum_{\pi \in \sortable{n}{t}}x^{\des (\pi)}
\end{align*}
be the  generating function for the descent statistic over $t$-stack-sortable permutations.

It is  known that $W_{n,1}(x)$ are the  Narayana polynomials
%~\cite{Stanley1989}
$$N_n(x) := \sum_{k=0}^{n-1}\frac{1}{n}\binom{n}{k}\binom{n}{k+1}x^k$$ and $W_{n,n-1}(x)$, being the set of all permutations of length $n$,  are the Eulerian polynomials $A_n(x)$.
% ~\cite{Kitaev2011Patterns}.
Jacquard and Schaeffer \cite{Jacquard1998bijective} gave the following formula in the case of $t=2$
\begin{align*}
	W_{n,2}(x)= \sum_{k=0}^{n-1}  \frac{(n+k)!(2n-k-1)!}{(k+1)!(n-k)!(2k+1)!(2n-2k-1)!} x^k,
\end{align*}
which also counts certain {\em planar maps} and so-called $\beta(0,1)$-trees encoding them \cite{Kitaev2011Patterns}.
For the case of $t=n-2$, Br{\"a}nd{\'e}n  \cite[Sec. 5]{Braenden2006linear} proved that
\begin{align}\label{n-2}
	W_{n,n-2}(x)=A_{n}(x)-x\, A_{n-2}(x).
\end{align}
%By using certain real-rootedness preserving linear operator, Br\"and\'en  proved the real-rootedness of $W_{n,n-2}(x)$.

Let $W_{n,t}^{p}(x)$ denote the generating function for descents over permutations in  $\sortable{n}{t}$ that avoid a pattern $p$. It is known \cite{BKLPRW} that there are
\begin{align}\label{1-st-sort-132}
	\frac{1}{k+1}{n-1\choose k}{n\choose k}
\end{align}
$p$-avoiding permutations of length $n$ with $k$ descents if $p\in\{132, 213, 231, 312\}$, which in particular gives the distribution of descents over 1-stack-sortable permutations (which are precisely 231-avoiding permutations). On the other hand, the distribution of descents over 123-avoiding permutations is given by the following formula \cite{BBS2010,BKLPRW}, where $t$ and $x$ correspond to the length and the number of descents,
\begin{align}
	1+ \sum_{n=1} t^n {\sum_{\substack{\pi\in\mathfrak S_n\\ \pi\text{ avoids }123}}x^{\des(\pi)}} = \frac{-1+2tx+2t^2x-2tx^2-4t^2x^2+2t^2x^3+\sqrt{1-4tx-4t^2x+4t^2x^2}}{2tx^2(tx-1-t)}. \notag
\end{align}
For $k\geq1$, the number of 321-avoiding permutations of length $n$ with $k$ descents \cite[A091156]{oeis} is
\begin{align}\label{321-des}
	\frac{1}{n+1}{n+1\choose k}\sum_{j=0}^{n-2k}{k+j-1\choose k-1}{n+1-k \choose n-2k-j}.
\end{align}
The number of descents plus the number of ascents in a permutation of length $n$ is $n-1$, and reading all 321-avoiding permutations in the reverse order gives all 123-avoiding permutations. Hence, we have the following distribution of descents over 123-avoiding permutations of length $n$,
\begin{align}\label{123-des}
	\frac{1}{n+1}{n+1\choose k+2}\sum_{j=0}^{2k-n+2}{n-k+j-2\choose j}{k+2 \choose n-k+j}.
\end{align}
A particular result of Bukata at el.~\cite{BKLPRW} states that the number of 213-avoiding 1-stack-sortable (i.e.\ 213-avoiding and 231-avoiding) permutations of length $n$ with $k$ descents is given by ${n-1\choose k}$. Egge and Mansour \cite{EggeMansour} studied 132-avoiding 2-stack-sortable permutations.

In this paper, we  give formulas for  $W_{n,n-3}(x)$ and $W_{n,n-4}(x)$ in terms of  Eulerian polynomials, which generalize the enumeration results of Claesson,  Dukes, and Steingr\'\i msson \cite{Claesson2010Permutations}. Moreover, we derive explicit formulas for  $W^{p}_{n,n-2}(x)$, $W^{p}_{n,n-3}(x)$ and $W^{p}_{n,n-4}(x)$, where $p$ is any permutation of length~3, which extend the studies in \cite{BBS2010,BKLPRW,EggeMansour}. Setting $x=1$ in these formulas gives us enumeration of $p$-avoiding $(n-2)$-, $(n-3)$- and $(n-4)$-stack-sortable permutations for $p\in\mathfrak{S}_3$ (see Table~\ref{length-3-pattern-avoidance}).

The paper is organized as follows. In Section~\ref{main-results-sec} we list the main enumerative results in the paper and introduce the stack-sorting complexity. In Sections~\ref{sect-2} and~\ref{sect-3}  we prove Theorems~\ref{main-thm1} and~\ref{main-thm2}, respectively. In Section~\ref{sec-remaining-proofs} we provide proofs of Theorems~\ref{231-thm}--\ref{312-thm}. Finally, in Section~\ref{final-sect} we give some concluding remarks and state two conjectures.

\section{The main results in the paper}\label{main-results-sec}

The main results in this paper are stated in the following theorems.
\begin{thm}\label{main-thm1}
	For $n\ge 4$, we have that
	\begin{align} \label{main-eq1}
		W_{n,n-3}(x)= & \,  A_n(x) - \frac{5}{2} x A_{n-2}(x)- \left( n-\frac{3}{2}\right) x(x+1)A_{n-3}(x).
	\end{align}
\end{thm}

\begin{thm}\label{main-thm2}
	For $n\ge 4$, we have that
	\begin{align}
		W_{n,n-4}(x) =
		\left\{\begin{array}{ll}
			1 & \mbox{if } n=4,\\
			x^4+10x^3+20x^2+10x+1 & \mbox{if } n=5,\\[2pt]
			{\begin{aligned}
			&\frac{1}{6}\Big[6A_n(x)-26xA_{n-2}(x)
			-(15n-30)x(x+1)A_{n-3}(x)\\
			&\qquad-\Big((3n^2-12n+10)(x^3+x)
			+(12n^2-48n+28)x^2\Big)A_{n-4}(x)\Big]
			\end{aligned}} & \mbox{if } n\geq6.
		\end{array}\right. \notag
	\end{align}
\end{thm}

\begin{thm}\label{231-thm}
	We have that
	\begin{align}
		W^{231}_{n,t}(x) = &
		\,  \displaystyle\sum_{k= 0}^{n-1}\frac{1}{k+1}{n-1\choose k}{n\choose k}x^k \quad \mbox{ for } 1 \le t\le n-1.	\notag
	\end{align}
\end{thm}

\begin{thm}\label{123-thm} For $n\ge 4$, we have that
	\begin{align}
		W^{123}_{n,n-2}(x) = &
		\,   W^{123}_{n,n-1}(x)  - x^{n-2}; \notag                   \\
		W^{123}_{n,n-3}(x) = &
		\,   W^{123}_{n,n-1}(x)  - (n-2)x^{n-3}-(n+1)x^{n-2}; \notag \\
		W^{123}_{n,n-4}(x) = &
		\,   W^{123}_{n,n-1}(x)  - 	\left\{\begin{array}{ll}
			                                   (2x+11x^2+x^3)                                              & \mbox{if }n=4     \\
			                                   (15x^2+16x^3)                                               & \mbox{if }n=5     \\
			                                   \left((n^2-n-3)x^{n-3} + \frac{n^2+n+2}{2}\, x^{n-2}\right) & \mbox{if }n\geq 6\end{array}\right. \notag
	\end{align}
	where $W^{123}_{n,n-1}(x)=\displaystyle\sum_{k\geq 0}\frac{1}{n+1}{n+1\choose k+2}\sum_{j=0}^{2k-n+2}{n-k+j-2\choose {j}}{k+2 \choose n-k+j}x^k$.
\end{thm}

\begin{thm}\label{321-thm} For $n\ge 4$, we have that
	\begin{align}
		W^{321}_{n,n-2}(x) = &
		\,   W^{321}_{n,n-1}(x)  - x; \notag                    \\
		W^{321}_{n,n-3}(x) = &
		\,   W^{321}_{n,n-1}(x)  \,  - (n+1)x -(n-3)x^2; \notag \\
		W^{321}_{n,n-4}(x) = &
		\,   W^{321}_{n,n-1}(x)  - \frac{n^2+n+2}{2}x-\left(n^2-n-10\right)x^2-\binom{n-4}{2}x^3 ; \notag
	\end{align}
	where $W^{321}_{n,n-1}(x)={1+\displaystyle\sum_{k\geq1}}\frac{1}{n+1}{n+1\choose k}\sum_{j=0}^{n-2k}{k+j-1\choose k-1}{n+1-k \choose n-2k-j}x^k$.
\end{thm}
Recall that the Narayana polynomials are defined as  $$N_{n}(x) = \displaystyle\sum_{k= 0}^{n-1}\frac{1}{k+1}{n-1\choose k}{n\choose k}x^k$$ for $n\ge 1$ and $N_{0}(x)=1$. It is known that
$$W^{213}_{n,n-1}(x)=W^{132}_{n,n-1}(x)= W^{312}_{n,n-1}(x)= N_n(x).$$
\begin{thm}\label{213-thm} For $n\ge 4$, we have that
	\begin{align}
		W^{213}_{n,n-2}(x) = &
		\,   N_{n}(x)  - x; \notag            \\
		W^{213}_{n,n-3}(x) = &
		\,   N_{n}(x)  - 3x-(2n-5)x^2; \notag \\
		W^{213}_{n,n-4}(x) = &
		\,   N_{n}(x)  - 	(6x+(8n-26)x^2+(2n-7)(n-3)x^3).  \notag
		%\left\{\begin{array}{ll}
		%(6x+6x^2+x^3) & \mbox{if }n=4 \\
		%(6x+(8n-26)x^2+(2n-7)(n-3)x^3) & \mbox{if }n\geq 5 \end{array}\right. \notag
	\end{align}
\end{thm}

\begin{thm}\label{132-thm} For $n\ge 4$, we have that
	\begin{align}
		W^{132}_{n,n-2}(x) = &
		\,   N_{n}(x)  - xN_{n-2}(x); \notag                     \\
		W^{132}_{n,n-3}(x) = &
		\,  N_{n}(x) - xN_{n-2}(x)  - 2(x+x^2)N_{n-3}(x); \notag \\
		W^{132}_{n,n-4}(x) = &
		\,   W^{132}_{n,n-3}(x)  -
		\left\{
		\begin{array}{ll}
			(3x+3x^2+x^3)            & \mbox{if }n=4     \\
			(3x+7x^2+3x^3)           & \mbox{if }n=5     \\
			(3x+8x^2+3x^3)N_{n-4}(x) & \mbox{if }n\geq 6
		\end{array}
		\right. \notag
	\end{align}
\end{thm}

\begin{thm}\label{312-thm} For $n\ge 4$, we have that
	\begin{align}
		W^{312}_{n,n-2}(x) = &
		\,   N_{n}(x)  - xN_{n-2}(x); \notag                  \\
		W^{312}_{n,n-3}(x) = &
		\,  N_{n}(x)  - 2xN_{n-2}(x)-x(x+1)N_{n-3}(x); \notag \\
		W^{312}_{n,n-4}(x) = &
		\,   W^{312}_{n,n-3}(x)  - 	\left\{\begin{array}{ll}
			                                   (3x+3x^2+x^3)                     & \mbox{if }n=4     \\
			                                   \big( xN_{n-2}(x)+ x^2 N_{n-3}(x) & \mbox{if }n\geq 5 \\  + (2x+2x^2+x^3) N_{n-4}(x)  \big)  &\end{array}\right. \notag
	\end{align}
\end{thm}

\begin{rem} Setting $x=1$ in Theorems~\ref{231-thm}--\ref{312-thm}, we obtain respective enumeration of $p$-avoiding $(n-2)$-, $(n-3)$- and $(n-4)$-stack-sortable permutations for $p\in\mathfrak{S}_3$ that is summarized in Table~\ref{length-3-pattern-avoidance}. \end{rem}

\begin{table}
	\begin{center}
		\begin{tabular}{|c|c|c|}
			\hline
			\hline
			$p$ & $(n-2)$-stack-sortable                                                                                                                                         & $(n-3)$-stack-sortable                                                                                                                                                                                   \\
			\hline
			\hline
			123 & \begin{tabular}{c} \footnotesize{13, 41, 131, 428, 1429, 4861, 16795, ... }\\ $C_n-1$; \cite[A001453]{oeis}  \end{tabular}                                     & \begin{tabular}{c}  \footnotesize{7, 33, 121, 416, 1415, 4845, 16777,  ...}\\  $C_n-2n+1$ \end{tabular}                                                                                                  \\[4mm]
			\hline
			132 & \begin{tabular}{c} \footnotesize{12, 37, 118, 387, 1298, 4433, 15366, ...}\\ $C_{n}-C_{n-2}$; \cite[A280891]{oeis} \end{tabular}  & \begin{tabular}{c} \footnotesize{8, 29, 98, 331, 1130, 3905, 13650, ...}\\ $C_{n} - C_{n-2}  - 4C_{n-3}$ \end{tabular} \\
			\hline
			213 & \begin{tabular}{c} \footnotesize{13, 41, 131, 428, 1429, 4861, 16795, ...} \\ $C_n-1$; \cite[A001453]{oeis}  \end{tabular}                                     & \begin{tabular}{c} \footnotesize{8, 34, 122, 417, 1416, 4846, 16778, ...} \\ $C_n-2(n-1)$ \end{tabular}                                                                                                  \\
			\hline
			231 & \begin{tabular}{c} \footnotesize{14, 42, 132, 429, 1430, 4862, 16796, ...} \\ $C_n$; \cite[A000108]{oeis}  \end{tabular}                                       & \begin{tabular}{c} \footnotesize{14, 42, 132, 429, 1430, 4862, 16796, ...} \\ $C_n$; \cite[A000108]{oeis}  \end{tabular}                                                                                 \\
			\hline
			312 & \begin{tabular}{c} \footnotesize{12, 37, 118, 387, 1298, 4433, 15366, ...} \\ $C_{n}-C_{n-2}$; \cite[A280891]{oeis} \end{tabular} & \begin{tabular}{c} \footnotesize{8, 28, 94, 317, 1082, 3740, 13078, ...} \\ $C_n-2(C_{n-2}+C_{n-3})$ \end{tabular}                                                            \\
			\hline
			321 & \begin{tabular}{c} \footnotesize{13, 41, 131, 428, 1429, 4861, 16795, ...} \\ $C_n-1$; \cite[A001453]{oeis} \end{tabular}                                      & \begin{tabular}{c} \footnotesize{8, 34, 122, 417, 1416, 4846, 16778, ...} \\ $C_n-2(n-1)$ \end{tabular}                                                                                                  \\
			\hline
		\end{tabular}

		\vspace{0.3cm}

		\begin{tabular}{|c|c|}
			\hline
			\hline
			$p$ & $(n-4)$-stack-sortable                                                                                                                                         \\
			\hline
			\hline
			123 & \begin{tabular}{c}  \footnotesize{0, 11, 83, 361, 1340, 4747, 16653, ... }\\ $C_n-\frac{3n^2-n-4}{2}$, $n\geq 6$ \end{tabular}                         \\[4mm]
			\hline
			132 & \begin{tabular}{c} \footnotesize{1, 16, 70, 261, 934, 3317, 11802, ...}\\ {$C_n-C_{n-2}-4C_{n-3}-14C_{n-4}$, $n\geq6$}\end{tabular}                                              \\
			\hline
			213 & \begin{tabular}{c} \footnotesize{1, 16, 89, 365, 1341, 4744, 16645, ...} \\ $C_n-(2n^2-5n+1)$ \end{tabular}                                                    \\
			\hline
			231 & \begin{tabular}{c} \footnotesize{1, 42, 132, 429, 1430, 4862, 16796, ...}\\ {$C_n$, $n\geq5$}; \cite[A000108]{oeis}   \end{tabular}                                        \\
			\hline
			312 & \begin{tabular}{c} \footnotesize{1, 16, 65, 236, 838, 2969, 10559, ...}\\ {$C_n-3C_{n-2}-3C_{n-3}-5C_{n-4}$, $n\geq5$} \end{tabular} \\
			\hline
			321 & \begin{tabular}{c} \footnotesize{1, 16, 89, 365, 1341, 4744, 16645, ...} \\ $C_n-(2n^2-5n+1)$  \end{tabular}                                                   \\
			\hline
		\end{tabular}

		\caption{Enumeration of $p$-avoiding $t$-stack-sortable permutations of length $n\geq 4$ for $t\in\{n-4,n-3,n-2\}$ where $C_n=\frac{1}{n+1}\binom{2n}{n}$ is the $n$-th Catalan number. For any $p$ of length 3, $(n-1)$-stack-sortable permutations are counted by $C_n$. Note that we have the same formulas for $p=213$ and $p=321$. In Section~\ref{final-sect} we state Conjecture~\ref{213-321-conj} generalizing this observation along with Conjecture~\ref{conj2} linking our permutations to the 321-machine considered in~\cite{Cerbai2020Stack}.}\label{length-3-pattern-avoidance}
	\end{center}
\end{table}

Following~\cite{Claesson2010Permutations}, let the
\emph{(stack-sorting) complexity} of $\pi$ be the
smallest integer $t$ such that $\pi$ is $t$-stack-sortable and
denote by  $\exactly{n}{t}$ (resp.,  $\exactly{n}{t}^p$) the set of permutations of length $n$ (resp., avoiding a pattern $p$) with  complexity $t$.
Denote  by
$$E^p_{n, t}(x)= \sum_{\pi \in\exactly{n}{t}^p}x^{\des (\pi)},$$
the  generating function for descents over $\exactly{n}{t}^p$, where $p$ can be omitted.
Then, we have that, for $2 \le t \le n-1$,
\begin{align}\label{eq-W-E}
	W^{p}_{n, t-1}(x)   =W^{p}_{n, t}(x)  - E^{p}_{n, t}(x)
\end{align}
where again $p$ can be omitted (if there is no pattern to avoid). Hence, together with~\eqref{n-2}, in order to compute $W_{n,n-3}(x)$, it suffices to find $E_{n,n-2}(x)$.
Similarly, once \eqref{main-eq1} is proved, the formula for $W_{n,n-4}(x)$ will be obtained from the formula for $E_{n,n-3}(x)$ given in Theorem~\ref{main}. The same approach is applied to compute $W^p_{n,n-2}(x)$, $W^p_{n,n-3}(x)$, $W^p_{n,n-4}(x)$ for $p\in\mathfrak{S}_3$ in Section~\ref{sec-remaining-proofs}.

\section{Proof of Theorem \ref{main-thm1}}\label{sect-2}

In this section, we shall prove Theorem~\ref{main-thm1}, which gives a formula for $W_{n,n-3}(x)$.
By \eqref{eq-W-E}, it suffices to compute $E_{n,n-2}(x)$.

Let us first review the structure of $\exactly{n}{n-2}$.
For convenience of  representation of certain sets of permutations over the alphabet $\{1,2,\ldots\}$, we use the  computer science notation of \emph{glob patterns}.  An asterisk ($\st$) stands
for any rearrangement of a subset of $[n]$ (including the empty permutation), and  a question mark
($\q$) stands for any single element in $[n]$.
For a word $w$ over the alphabet $\{\st,\q\}\cup\{1,2,\dots\}$, $\glob{w}$ denotes the set of
permutations obtained from $w$ by all possible replacements of $\st$'s and $\q$'s.  Also,  let
$$\glob{w_1,\dots,w_k} := \glob{w_1}\cup\dots\cup\glob{w_k}.$$
West  \cite{West1990Permutations}  characterized $\exactly{n}{n-2}$ for all $n\geq 4$ as follows,
\begin{align}\label{West-n-2-exact}
	\exactly{n}{n-2} = \sym_n \cap
	\glob{\;\st n2,\;\st(n-1)1n,\; \st n1\q,\; \st n\q 1,
		\; \st n\st (n-2)\st (n-1)1\;}.
\end{align}

For convenience, we let $A_n \glob{w}$ be the  generating function for distribution of descents over  permutations of $[n]$ in $\glob{w}$.
Easy observations lead to the following equations, where the factor of $x$ (resp., $x^2$) records the extra descent (resp., two descents),
\begin{align*}
	A_n \glob{\st n2}     & = x A_{n-2}(x),         \\[-2pt]
	A_n \glob{\st(n-1)1n} & =  x A_{n-3}(x),        \\[-2pt]
	A_n \glob{\st n1\q}   & = (n-2) x A_{n-3}(x),   \\[-2pt]
	A_n \glob{\st n\q 1}  & = (n-2) x^2 A_{n-3}(x).
\end{align*}

We proceed by computing  the remaining case of $A_n \glob{\st n\st (n-2)\st (n-1)1}$.
\begin{lem}\label{lem-1}
	For any $n\ge 2$, we have  that
	$$A_n \glob{\st n\st (n-1)\st }
		= \frac{1}{2}\big( A_{n}(x) + (x-1) A_{n-1}(x) \big).$$
\end{lem}

\begin{proof}
	We note that the set of all permutations of $[n]$ is the disjoint union of the following four  sets
	\begin{align*}
		 & \ \glob{\st n \, ? *  (n-1)\st}   \bigcup \,    \glob{\st n  (n-1)\st} \bigcup \,  \glob{\st (n-1) \, ? * n \st}  \bigcup \,   \glob{\st (n-1) n \st}.
	\end{align*}
	%	\begin{align*}
	%	& \ \glob{\st n \, ? *  (n-1)\st,  \,  \st n  (n-1)\st, \,   \st (n-1) \, ? * n \st,   \st (n-1) n \st},
	%	\end{align*}
	Since
	\begin{align*}
		A_n \glob{\st n  (n-1)\st} =  x A_{n-1}(x), \quad
		A_n \glob{\st  (n-1) n \st} =   A_{n-1}(x), \\[1pt]
		A_n \glob{\st \, n \, ? *  (n-1)\st}   = A_n \glob{\st (n-1) \, ? * n \st},
	\end{align*}
	it follows that
	\begin{align*}
		2 A_n \glob{\st\, n \, ? *  (n-1)\st}  +  x A_{n-1}(x) + A_{n-1}(x)  = A_{n}(x),
	\end{align*}
	and thus
	\begin{align*}
		A_n \glob{\st\, n \, ? *  (n-1)\st}   = \frac{1}{2} \big(A_{n}(x)-(x+1)A_{n-1}(x)\big).
	\end{align*}
	From
	\begin{align*}
		\glob{\st\,  n *  (n-1)\st}    =  \glob{\st \,  n \, ? \, *  (n-1)\st} \bigcup \,  \glob{\st \, n  (n-1)\, \st}
	\end{align*}
	we have that
	\begin{align*}
		A_n \glob{\st\,  n \, *  (n-1)\st} & =  A_n \glob{\st \, n \, ? \, *  (n-1)\st}  + A_n \glob{\st\,  n  (n-1)\, \st} \\[-5pt]
		                                   & = \cfrac{1}{2}  \left( A_{n}(x) + (x-1) A_{n-1}(x) \right).
	\end{align*}
	This completes the proof.
\end{proof}

Now we are able to obtain an expression for $ E_{n,n-2}(x)$.
\begin{thm}\label{thm-1}
	For all $n \ge 4$,
	\begin{align} \label{e:n-2}
		E_{n,n-2}(x)= & \frac{1}{2} x \big( 3 A_{n-2}(x) + (2 n-3) (x+1) A_{n-3}(x) \big).
	\end{align}
\end{thm}

\begin{proof}
	Note that  for any permutation in $\sym_n \cap    \glob{\st n\st (n-2)\st (n-1)1}$, the leftmost $n-2$ elements have the same pattern as any permutation in
	$\sym_{n-2} \cap    \glob{\st (n-2) \st (n-3)\st}$.
	It then follows from  Lemma \ref{lem-1} that
	\begin{align*}
		A_n \glob{\st n\st (n-2)\st (n-1)1} & = \frac{1}{2} x \left( A_{n-2}(x) + (x-1) A_{n-3}(x) \right)
	\end{align*}
	where the factor of $x$ corresponds to the descent $(n-1)1$. Summing up the polynomials corresponding to the five words in \eqref{West-n-2-exact}, we obtain the desired result.
	This completes the proof.
\end{proof}

Since $W_{n,n-3}(x)   = W_{n,n-2}(x) - E_{n,n-2}(x)$,
formula \eqref{main-eq1} follows immediately from \eqref{n-2} and \eqref{e:n-2}, which completes the proof of  Theorem \ref{main-thm1}.

\section{Proof of Theorem \ref{main-thm2}}\label{sect-3}

In this section, we shall prove Theorem \ref{main-thm2}, which gives a formula for $W_{n,n-4}(x)$.
The cases $n=4,5$ are verified directly. Hence, throughout the classification argument below, we assume $n\geq6$.
By \eqref{eq-W-E}, it suffices to compute $E_{n,n-3}(x)$.

Our proof is based on the classification of permutations in $\exactly{n}{n-3}$,
listed in  Tables~\ref{thelist1} and~\ref{thelist2},  that were given by Claesson,  Dukes, and
Steingr\'\i msson~\cite{Claesson2010Permutations}. We shall compute the descent generating polynomial $A_n \glob{w}$  for each type of $w$. To do this, we divide all the cases into {\em four} classes.

\begin{table}[h!]
	%{		\renewcommand{\arraystretch}{1.05}
	$\begin{array}{|l|c|l|} \hline\hline
			                                                   &                                                                    &                                                            \\[-2ex]
			\mbox{Case}                                        & \mbox{Type of } w                                                  & \mbox{Descent gen. polynomial }  A_n \glob{w}              \\[-2ex]
			                                                   &                                                                    &                                                            \\\hline\hline
			\multicolumn{3}{l}{}                                                                                                                                                                 \\[-3ex]
			\multicolumn{3}{l}{\pi_n=n}                                                                                                                                                          \\ \hline
			\mbox{1(a)}                                        & \st (n-1)2n                                                        & x A_{n-3}(x)                                               \\  \hline
			\mbox{1(b)}                                        & \st (n-1)\q 1n                                                     & (n-3) x^2 A_{n-4}(x)                                       \\   \hline
			\mbox{1(c)}                                        & \st (n-1)1\q n                                                     & (n-3) x A_{n-4}(x)                                         \\  \hline
			\mbox{1(d)}                                        & \st (n-1) \st (n-3)\st (n-2) 1 n                                   & \frac{1}{2} x \left( A_{n-3}(x) + (x-1) A_{n-4}(x) \right) \\[3pt] \hline
			\mbox{1(e)}                                        & \st (n-2)1(n-1)n                                                   & x A_{n-4}(x)                                               \\ \hline
			\multicolumn{3}{l}{}                                                                                                                                                                 \\[-3ex]
			\multicolumn{3}{l}{\pi_{n-1}=n}                                                                                                                                                      \\ \hline
			\mbox{2(a)}                                        & \st n3                                                             & x A_{n-2}(x)                                               \\  \hline
			\mbox{2(b)}                                        &
			\st (n-1)1n\q, \text{where}\ \q = 4, 5,\ldots, n-2 & (n-5) x^2 A_{n-4}(x)                                                                                                            \\  \hline
			\mbox{2(c)}                                        &
			\st (n-2)1n(n-1)                                   & x^2 A_{n-4}(x)                                                                                                                  \\ \hline
			\multicolumn{3}{l}{}                                                                                                                                                                 \\[-3ex]
			\multicolumn{3}{l}{\pi_{n-2}=n}                                                                                                                                                      \\ \hline
			\mbox{3(a)}                                        &
			\st n2\q, \text{ where}\   \q \neq 1               & (n-3)x  A_{n-3}(x)                                                                                                              \\ \hline
			\mbox{3(b)}                                        &
			\st n\q 2,  \text{ where}\  \q \neq 1              & (n-3)x^2 A_{n-3}(x)                                                                                                             \\[3pt] \hline
			\multicolumn{3}{l}{}                                                                                                                                                                 \\[-3ex]
			\multicolumn{3}{l}{\pi_{n-3}=n}                                                                                                                                                      \\ \hline
			\mbox{4(a)}                                        &
			\st n\q \q 1, \mbox{ but not } \st n(n-2)(n-1)1    & \left( \binom{n-2}{2}x^3 +(\binom{n-2}{2}-1)x^2 \right) A_{n-4}(x)                                                              \\[3pt] \hline
			\mbox{4(b)}                                        & \st n \q 1 \q                                                      & (n-2)(n-3) x^2 A_{n-4}(x)                                  \\ \hline
			\mbox{4(c)}                                        & \st n 1\q \q                                                       & \binom{n-2}{2}(x+x^2)A_{n-4}(x)                            \\[3pt] \hline
			\mbox{4(d)}                                        & \st n(n-2)(n-1)2                                                   & x^2 A_{n-4}(x)                                             \\ \hline
		\end{array}$\bigskip
	\caption{Permutations in $\exactly{n}{n-3}$ for {$n\geq6$}}
	\label{thelist1}
\end{table}

\begin{table}[h!]
	{		\renewcommand{\arraystretch}{1}
		$\begin{array}{|l|c|l|} \hline\hline
				                                      &                                                                                                &                                                                                 \\[-2ex]
				\mbox{Case}                           & \mbox{Type of } w                                                                              & \mbox{Descent generating polynomial }  A_n \glob{w}                             \\[-2ex]
				                                      &                                                                                                &                                                                                 \\\hline\hline

				\multicolumn{3}{l}{}                                                                                                                                                                                                     \\[-3ex]
				\multicolumn{3}{l}{\pi_{n-i}=n \mbox{ and }i> 3}                                                                                                                                                                         \\ \hline

				\mbox{5(a)}
				                                      & \st n A (n-2) B(n-1)2 , \mbox{ where $A\cup B\neq \emptyset$}
				                                      & \mbox{ \small  $ \frac{1}{2} x \left( A_{n-2}(x) + (x-1) A_{n-3}(x) \right) - x^2 A_{n-4}(x)$}                                                                                   \\[1.6ex] \hline

				\mbox{5(b)}                           & \st n\st (n-3)\st (n-1)(n-2)1                                                                  & \mbox{ \small  $\frac{1}{2} x^2  \left( A_{n-3}(x) + (x-1) A_{n-4}(x) \right)$} \\[1.6ex]  \hline

				\mbox{5(c)}                           &
				\st n \st (n-1) \st (n-3) \st (n-2)1  &
				\begin{array}{l}
					\mbox{ \small  $\frac{1}{6} x \left( A_{n-2}(x) + 3(x-1) A_{n-3}(x) +\right.$} \\[0.5ex]
					\mbox{ \small  $\left.  + 2(x-1)^2 A_{n-4}(x)  \right)$}
				\end{array}                                                            \\[2.6ex]	 \hline
				%	\mbox{ \small  $\frac{1}{6} x \left( A_{n-2}(x) + 3(x-1) A_{n-3}(x) + 2(x-1)^2 A_{n-4}(x) \right)$}
				% \\[2.6ex]

				\mbox{5(d)}                           & \st (n-1) \st n \st\!\left\{\!\!\!
				\begin{array}{c}
					(n-3) \st (n-4) \\
					(n-4) \st (n-3)
				\end{array}
				\!\!\!\right\}\!\st (n-2)1            &
				\mbox{\small {$\displaystyle
				\frac{x}{12}\left(A_{n-2}(x)-(x-1)^2A_{n-4}(x)\right)$}}                                                              \\[3.6ex]  \hline

				\mbox{5(e)}                           & \st n \st (n-2) \st (n-1)\!\left\{\!\!\!
				\begin{array}{c}
					1\q \\ \q 1
				\end{array}\!\!\!\right\}             & \mbox{ \small  $ (n-4)(\frac{x}{2}+\frac{x^2}{2})\left( A_{n-3}(x)+(x-1)A_{n-4}(x)\right)$}                                                                                      \\[2ex]  \hline

				\mbox{5(f)}                           & \st n\st (n-3)\st (n-1)1(n-2)                                                                  & \mbox{ \small  $ \frac{1}{2} x \left( A_{n-3}(x)+(x-1)A_{n-4}(x)\right)$}       \\[1.6ex]  \hline

				\mbox{5(g)}                           & \st n \st (n-3) \st (n-2)1(n-1)                                                                & \mbox{ \small  $\frac{1}{2}  x    \big( A_{n-3}(x)+(x-1)A_{n-4}(x)\big)$}       \\[1.6ex]   \hline

				\mbox{5(h)}                           & \st (n-2)\st n\st\!\left\{\!\!\!
				\begin{array}{c}
					(n-3)\st (n-4) \\
					(n-4)\st (n-3)
				\end{array}\!\!\!\right\}\!\st (n-1)1 &
				\mbox{\small {$\displaystyle
				\frac{x}{12}\left(A_{n-2}(x)-(x-1)^2A_{n-4}(x)\right)$}}                                                              \\ \hline
			\end{array}$\bigskip}
	\caption{Permutations in $\exactly{n}{n-3}$ for {$n\geq6$}}\label{thelist2}
\end{table}

\noindent \textbf{Class 1:} All cases in Table~\ref{thelist1}  except  case 1(d)

The descent generating functions for the cases in Table~\ref{thelist1}  except 1(d) can be obtained directly.
For example, consider type 4(a), namely $$\st n\q \q 1, \mbox{ but not } \st n(n-2)(n-1)1.$$
Since the letters $\q \q$ can be  chosen arbitrarily from $\{2, 3, \ldots, n-1 \}$,  and ordering them in decreasing order results in an extra descent, we obtain that
$$A_n \glob{ \st n\q \q 1 } = \left( \binom{n-2}{2}x^3 + \binom{n-2}{2}  x^2 \right) A_{n-4}(x), $$
and hence the descent generating polynomial for case 4(a) is
$$ \left( \binom{n-2}{2}x^3 +\left(\binom{n-2}{2}-1\right)x^2 \right) A_{n-4}(x). $$

\noindent \textbf{Class 2:} 1(d), 5(a), 5(b), 5(e), 5(f), and 5(g)

The  descent generating functions for cases 1(d), 5(a), 5(b), 5(e), 5(f), and 5(g) can be derived from Lemma \ref{lem-1}.
For instance, in the case of 5(a),  the type
$\st n A (n-2) B(n-1)2$ , where $A\cup B\neq \emptyset$, can be translated into the  form $$\st n\st (n-2) \st (n-1)2, \quad \mbox{  but not  } \quad   \st n (n-2)  (n-1)2,$$ and hence we obtain the  descent generating polynomial as desired.

\noindent \textbf{Class 3:} 5(c)

To deal with case 5(c), we need the following lemma.
\begin{lem}\label{lem-length-3}
	For all $n\ge 3$,
	\begin{align*}
		A_n \glob{\st n\st (n-1)\st (n-2) \st }
		= \frac{1}{6} \left( A_{n}(x) + 3 (x-1) A_{n-1}(x) + 2(x-1)^2 A_{n-2}(x) \right).
	\end{align*}
\end{lem}
\begin{proof}
	We note that the set of all permutations of $[n]$ is the disjoint union of the following three types, where $abc$ forms a permutation of $\{n-2,n-1,n\}$:
	
	(i) $\glob{\st\,a\,?\st\,b\,?\st\,c\,\st}$;
	
	(ii) $\glob{\st\,a\,b\,?\st\,c\,\st}$ or $\glob{\st\,a\,?\st\,b\,c\,\st}$;
	
	(iii) $\glob{\st\,a\,b\,c\,\st}$.
	
	The strict-gap form of the computation in Lemma~\ref{lem-1} is
	\[
	A_{n-1}\glob{\st(n-1)\,?\st(n-2)\st}
	=\frac12\left(A_{n-1}(x)-(x+1)A_{n-2}(x)\right).
	\]
	This polynomial is unchanged if the two distinguished top letters are interchanged, because the enforced gap makes all their adjacent neighbours smaller.
	Thus
	\begin{align*}
		A(\sym_n\cap\glob{\st a\,b\,?\st\,c\,\st})
		&=x^{\chi(a>b)}\frac12\left(A_{n-1}(x)-(x+1)A_{n-2}(x)\right),\\
		A(\sym_n\cap\glob{\st a\,?\st\,b\,c\,\st})
		&=x^{\chi(b>c)}\frac12\left(A_{n-1}(x)-(x+1)A_{n-2}(x)\right),\\
		A(\sym_n\cap\glob{\st a\,b\,c\,\st})
		&=x^{\chi(a>b)+\chi(b>c)}A_{n-2}(x),
	\end{align*}
	where $\chi(\cdot)$ is $1$ if the statement is true and $0$ otherwise. Since relabelling the three separated distinguished letters in type (i) does not change descents, the six orders have the same generating polynomial. Hence
	\begin{align*}
		&6A(\sym_n\cap\glob{\st a\,?\st\,b\,?\st\,c\,\st})
		+6(x+1)\frac12\left(A_{n-1}(x)-(x+1)A_{n-2}(x)\right)\\
		&\qquad +(x^2+4x+1)A_{n-2}(x)=A_n(x),
	\end{align*}
	and therefore
	\begin{align*}
		A_n\glob{\st a\,?\st\,b\,?\st\,c\,\st}
		=\frac16\left(A_n(x)-3(x+1)A_{n-1}(x)+2(x^2+x+1)A_{n-2}(x)\right).
	\end{align*}
	
	The set counted in the lemma is the disjoint union
	\begin{align*}
		&\glob{\st n\,?\st(n-1)\,?\st(n-2)\st}
		\mathbin{\bigcup}\glob{\st n(n-1)\,?\st(n-2)\st}\\
		&\qquad\mathbin{\bigcup}\glob{\st n\,?\st(n-1)(n-2)\st}
		\mathbin{\bigcup}\glob{\st n(n-1)(n-2)\st}.
	\end{align*}
	Consequently,
	\begin{align*}
		A_n\glob{\st n\st(n-1)\st(n-2)\st}
		={}&\frac16\left(A_n(x)-3(x+1)A_{n-1}(x)+2(x^2+x+1)A_{n-2}(x)\right)\\
		&+2\frac{x}{2}\left(A_{n-1}(x)-(x+1)A_{n-2}(x)\right)+x^2A_{n-2}(x)\\
		={}&\frac16\left(A_n(x)+3(x-1)A_{n-1}(x)+2(x-1)^2A_{n-2}(x)\right).
	\end{align*}
	This completes the proof.
\end{proof}

\noindent \textbf{Class 4:} 5(d) and 5(h)

\begingroup
Let $R_n(x)$ denote the descent generating polynomial of either case 5(d) or case 5(h). We first consider 5(d), and write its permutations as
\[
B_0\,(n-1)\,B_1\,n\,B_2\,a\,B_3\,b\,B_4\,(n-2)1,
\qquad \{a,b\}=\{n-3,n-4\},
\]
where $(B_0,\ldots,B_4)$ is an ordered partition of the remaining labels. As in the proof of Lemma~\ref{lem-length-3}, an ordinary block contributes $\mathcal A(z)$, whereas a block that creates an additional boundary descent precisely when it is nonempty contributes
\[
Q(z)=1+x\bigl(\mathcal A(z)-1\bigr)=1-x+x\mathcal A(z).
\]
If $(a,b)=(n-3,n-4)$, the three forced descents are the ones beginning at $n$, $n-3$, and $n-2$; the blocks $B_1$ and $B_4$ are of $Q$-type, while $B_0,B_2,B_3$ are ordinary. This gives $x^3\mathcal A^3Q^2$. If $(a,b)=(n-4,n-3)$, only the descents beginning at $n$ and $n-2$ are forced, and $B_1,B_3,B_4$ are of $Q$-type, giving $x^2\mathcal A^2Q^3$. The same block analysis applies to 5(h), with $n-2$ and $n-1$ interchanged in the two fixed outer positions. Hence
\[
\sum_{m\geq0}R_{m+6}(x)\frac{z^m}{m!}
=x^3\mathcal A(z)^3Q(z)^2+x^2\mathcal A(z)^2Q(z)^3.
\]
A direct differentiation using $\mathcal A'=\mathcal A Q$ yields
\[
x^3\mathcal A^3Q^2+x^2\mathcal A^2Q^3
=\frac{x}{12}\left(\mathcal A^{(4)}-(x-1)^2\mathcal A''\right).
\]
Therefore, for each of 5(d) and 5(h),
\[
R_n(x)=\frac{x}{12}\left(A_{n-2}(x)-(x-1)^2A_{n-4}(x)\right).
\]
\endgroup
%consider the descent generating polynoimal of the type
%$\st a ? \st b ? \st c ? \st d \st$, where $a,b,c,d$ is any permutation of $\{n-3,n-2,n-1,n\}$.

Now we complete the proof of the formulas of the descent generating polynomials of all types listed in  Tables~\ref{thelist1} and~\ref{thelist2}.
By summing  up all these polynomials, we obtain the following expression for $E_{n,n-3}(x)$.

\begin{thm}\label{main}
	For all $n\geq 6$, the descent generating polynomial of permutations in $\exactly{n}{n-3}$ is
	\begingroup
	\begin{align}
		E_{n,n-3}(x)=\frac16\Big[&11xA_{n-2}(x)
		+(9n-21)x(x+1)A_{n-3}(x)\notag\\
		&+\Big((3n^2-12n+10)(x^3+x)
		+(12n^2-48n+28)x^2\Big)A_{n-4}(x)\Big]. \notag
	\end{align}
	\endgroup
\end{thm}
Since $W_{n,n-4}(x)  = W_{n,n-3}(x) - E_{n,n-3}(x)$, Theorem~\ref{main-thm2} follows immediately from \eqref{main-eq1} and Theorem~\ref{main}.
%which completes the proof of Theorem~\ref{main-thm2}.

\section{Proofs of Theorems~\ref{231-thm}--\ref{312-thm}}\label{sec-remaining-proofs}

We note that Theorem~\ref{231-thm} is an immediate corollary to formula (\ref{1-st-sort-132}) and the observation that  1-stack-sortable permutations are precisely 231-avoiding permutations.

For the other theorems, we can follow the same steps as in the proofs of Theorems~\ref{main-thm1} and~\ref{main-thm2}, and apply suitable known formulas for $W^{p}_{n,n-1}$ in question. Namely, we will iterate (\ref{eq-W-E}) for $t=n-1, n-2, n-3$ for each pattern $p$ in question. Additionally, we need to make sure that the forbidden pattern does not occur in them. Below, we provide brief (omitting justifications for easier parts) proofs of Theorems~\ref{123-thm}--\ref{312-thm}, where $\downarrow$ and $\uparrow$ stand for decreasing and increasing permutations, respectively.

The boundary cases not covered by the classification tables are verified directly.

\subsection{Proof of Theorem~\ref{123-thm}}

\begin{table}[h!]
	{		\renewcommand{\arraystretch}{1.1}
		$\begin{array}{|l|c|l|} \hline\hline
				                                             &                                                                &                                                                \\[-2ex]
				\mbox{Case}                                  & \mbox{Type of } w                                              & \mbox{Descent generating polynomial }  A_n \glob{w}            \\[-2ex]
				                                             &                                                                &                                                                \\\hline\hline
				%	\mbox{1(b)} & 3214 & x^2 \mbox{ (for } n=4 \mbox{ only) }\\  \hline
				\mbox{2(a)}                                  & \downarrow n3                                                  & x^{n-2}                                                        \\  \hline
				%	\mbox{2(c)} &
				%	2143 &  x^2  \mbox{ (for } n=4 \mbox{ only) }\\ \hline
				%	\mbox{3(a)} &
				%	\downarrow n21 & x^{n-2}\ \mbox{($\q$ must be 1)}\\ \hline
				\mbox{3(b)}                                  &
				\downarrow n\q 2,  \text{ where}\  \q \neq 1 & (n-3)x^{n-2}                                                                                                                    \\  \hline
				\mbox{4(a)}                                  &
				{\downarrow n\q \q 1}                 &
				%	 x^3  \mbox{ if } n=4;  \mbox{ otherwise }
				\frac{(n-2)(n-3)}{2}x^{n-2}+ (n-3)x^{n-3} \                                                                                                                                    \\ \hline
				\mbox{4(b)}                                  & \downarrow n \q 1 \q                                           &
				%	2x^2 \mbox{ if } n=4; \mbox{ otherwise }
				\frac{n(n-3)}{2}x^{n-3}
				\                                                                                                                                                                              \\
				                                             &                                                                &
				\mbox{(consider the cases }2\in\downarrow \mbox{ and } 2\not\in\downarrow)                                                                                                     \\ \hline
				\mbox{4(c)}                                  & \downarrow n 1ba,\ \text{ where}\ \ b>a                        &
				%	x^2 \mbox{ if } n=4; \mbox{ otherwise }
				\frac{(n-2)(n-3)}{2} x^{n-3}                                                                                                                                                   \\ \hline
				\mbox{5(a)}
				                                             & n (n-2) \downarrow(n-1)2 , \mbox{ $\downarrow \neq \emptyset$}
				                                             & x^{n-2}                                                                                                                         \\ \hline

				\mbox{5(b)}                                  & n(n-3)\downarrow (n-1)(n-2)1                                   & x^{n-2}                                                        \\  \hline

				\mbox{5(c)}                                  &
				n (n-1) (n-3) \downarrow (n-2)1              & x^{n-2}
				\\	 \hline
				%	\mbox{ \small  $\frac{1}{6} x \left( A_{n-2}(x) + 3(x-1) A_{n-3}(x) + 2(x-1)^2 A_{n-4}(x) \right)$}
				% \\[2.6ex]

				\mbox{5(d)}                                  & (n-1)n(n-3)(n-4)\downarrow (n-2)1                              &
				x^{n-3}                                                                                                                                                                        \\  \hline

				\mbox{5(e)}                                  & n (n-2) \downarrow (n-1)\!\left\{\!\!\!
				\begin{array}{c}
					1\q \\ \q 1
				\end{array}\!\!\!\right\}                    & (n-4)(x^{n-3}+x^{n-2})                                                                                                          \\  \hline

				\mbox{5(f)}                                  & n(n-3)\downarrow (n-1)1(n-2)                                   & x^{n-3}                                                        \\  \hline
				\mbox{5(h)}                                  & (n-2)n(n-3)(n-4)\downarrow (n-1)1                              &
				x^{n-3}                                                                                                                                                                        \\   \hline
				%	\mbox{ \small  $\frac{1}{6} x \left( A_{n-2}(x) + 3(x-1) A_{n-3}(x) + 2(x-1)^2 A_{n-4}(x) \right)$}
				% \\[2.6ex]
			\end{array}$\bigskip}
	\caption{Permutations in $\exactly{n}{n-3}^{123}$ for $n\geq 6$}\label{thelist-123}
\end{table}

Considering (\ref{West-n-2-exact}) we can see that for $n\geq 4$,
\begin{align}\label{West-n-2-exact-123}
	\exactly{n}{n-2}^{123} = \sym_n \cap
	\glob{\;\downarrow n2,\; \downarrow n1\q,\; \downarrow n\q 1,
		\; n (n-2)\downarrow (n-1)1\;}
\end{align}
where, for example, $\st(n-1)1n$ disappears in the unrestricted set $\exactly{n}{n-2}$ because $n\geq 4$ and any element in the $\st$ along with $(n-1)$ and $n$ would form the forbidden pattern 123. Also, because of $(n-1)$,  the leftmost two $\st$'s in  $\st n\st (n-2)\st (n-1)1$ in $\exactly{n}{n-2}$ must be empty, while the remaining $\st$ must be a decreasing permutation. And so on.

It is easy to see that $A_n \glob{\downarrow n2}=x^{n-2}$, $A_n \glob{\downarrow n1\q}=(n-2)x^{n-3}$, {$A_n \glob{\downarrow n\q1}=(n-2)x^{n-2}$} and $A_n \glob{n (n-2)\downarrow (n-1)1}=x^{n-2}$, and hence from (\ref{West-n-2-exact-123}) we have
$$E^{123}_{n,n-2}(x)=(n-2)x^{n-3}+nx^{n-2}.$$
Now, since $\exactly{n}{n-1}=\sym_n \cap \glob{\st n1}$ (\cite{West1990Permutations}),  $\exactly{n}{n-1}^{123}=\sym_n \cap \glob{\downarrow n1}$ and $E^{123}_{n,n-1}(x)=x^{n-2}$. Since $\sortable{n}{n-1}=\sym_n$, formula (\ref{123-des}) can be used to find $W^{123}_{n,n-1}(x)$, and then iteration of (\ref{eq-W-E}) for $t=n-1, n-2$ gives the desired formulas for $W^{123}_{n,n-2}(x)$ and $W^{123}_{n,n-3}(x)$.

To complete the proof of Theorem~\ref{123-thm}, we need to compute $E^{123}_{n,n-3}(x)$ and apply (\ref{eq-W-E}) for $t=n-3$. We analyse Tables~\ref{thelist1} and~\ref{thelist2} presenting $\exactly{n}{n-3}$. We see that to avoid the pattern 123, the following cases are impossible when $n\ge 6$: 1(a), 1(b),  1(c), 1(d), 1(e), 2(b) (there is no place for 3), 2(c), 3(a), 4(d) (there is no place for 1 even for $n=4$), 5(g). The remaining cases are listed in Table~\ref{thelist-123}, along with several brief comments, and they give the desired formula for $W^{123}_{n,n-4}(x)$. In case 4(a), the excluded subclass $\st n(n-2)(n-1)1$ cannot avoid 123, so the type reduces to $\downarrow n\q\q1$.

\subsection{Proof of Theorem~\ref{321-thm}}

Considering (\ref{West-n-2-exact}) we can see that for $n\geq 4$,
\begin{align}\label{West-n-2-exact-321}
	\exactly{n}{n-2}^{321} = \sym_n \cap
	\glob{\;1\uparrow n2,\; An2,\; \uparrow (n-1)1n,\; \uparrow n1\q\;}
\end{align}
where, the element 1 in $A$ is not in the leftmost position, while all other elements in $A$ are increasing (because of the element 2). Indeed, for example, in the unrestricted set $\exactly{n}{n-2}$ $\st(n-1)1n$ becomes $\uparrow(n-1)1n$  because of the element 1, and $\st n\q1$ disappears  because no matter what $\q$ is, $n\q1$ is an occurrence of the pattern 321.

It is easy to see that $A_n \glob{1\uparrow n2}=x$, $A_n \glob{An2}=(n-3)x^2$, $A_n \glob{\uparrow (n-1)1n}=x$ and $A_n \glob{\uparrow n1\q}=(n-2)x$, and hence from (\ref{West-n-2-exact-321}) we have
$$E^{321}_{n,n-2}(x)=nx+(n-3)x^2.$$
Now, since $\exactly{n}{n-1}=\sym_n \cap \glob{\st n1}$ (\cite{West1990Permutations}),  $\exactly{n}{n-1}^{321}=\sym_n \cap \glob{\uparrow n1}$ and $E^{321}_{n,n-1}(x)=x$. Since $\sortable{n}{n-1}=\sym_n$, formula (\ref{321-des}) can be used to find $W^{321}_{n,n-1}(x)$, and then iteration of (\ref{eq-W-E}) for $t=n-1, n-2$ gives the desired formulas for $W^{321}_{n,n-2}(x)$ and $W^{321}_{n,n-3}(x)$.

To complete the proof of Theorem~\ref{321-thm}, we need to compute $E^{321}_{n,n-3}(x)$ and apply (\ref{eq-W-E}) for $t=n-3$. We analyse Tables~\ref{thelist1} and~\ref{thelist2} presenting $\exactly{n}{n-3}$. We see that to avoid the pattern 321, none of the cases in Table~\ref{thelist2} are possible, and the following cases are impossible too: 1(b), 1(d), 3(b), 4(a), 4(b) and 4(d). The remaining cases are listed in Table~\ref{thelist-321} and they give the desired formula for $W^{321}_{n,n-4}(x)$. Here are our comments for three (more involved) cases in Table~\ref{thelist-321}.
\begin{itemize}
	\item[1(a)] Because of the element 2, all elements but 1 in $\st$ must increase. The case of $1\uparrow(n-1)2n$ gives the term of $x$, and placing 1 differently in $\st$ gives the term of $(n-4)x^2$.
	\item[2(a)] Because of the element 3, all elements  greater than 3  in $\st$ must be increasing. Note that if $\st$ begins with 2, then we can place 1 in $(n-3)$ places that corresponds to the term of $(n-3)x^2$. Otherwise, to avoid the pattern 321, 1 must precede 2. We distinguish four possible cases here: $12\uparrow n3$ giving $x$; $\st 12 \st n3$, where to the left of 12 there is at least one element that gives $(n-4)x^2$; $1\st 2\st n3$, where there is at least one element between 1 and 2, that gives $(n-4)x^2$, and the remaining case of 1 preceding 2, 1 being not the leftmost element, and 1 and 2 staying not together, which gives the term of ${n-4 \choose 2}x^3$.
	\item[3(a)] Because of the element $2$, all elements in $\st$ that are  greater than 3  must be increasing and $\q\neq 1$, because otherwise $n21$ is an occurrence of the pattern 321. The case of $1\uparrow n2\q$ gives the term of $(n-3)x$ and the remaining cases of placing 1 differently and choosing $\q$ give the term of {$(n-3)(n-4)x^2$}.
\end{itemize}

\begin{table}[h!]
	%{		\renewcommand{\arraystretch}{1.05}
	$\begin{array}{|l|c|l|} \hline\hline
			            &                                                             &                                         \\[-2ex]
			\mbox{Case} & \mbox{Type of } w                                           & \mbox{Descent gen. pol. }  A_n \glob{w} \\[-2ex]
			            &                                                             &                                         \\\hline\hline
			\mbox{1(a)} & \st (n-1)2n                                                 & x + (n-4)x^2                            \\  \hline
			\mbox{1(c)} & \uparrow (n-1)1\q n                                         & (n-3)x                                  \\  \hline
			\mbox{1(e)} & \uparrow (n-2)1(n-1)n                                       & x                                       \\ \hline
			\mbox{2(a)} & \st n3                                                      & x + (3n-11)x^2 + {n-4 \choose 2}x^3     \\ \hline
			\mbox{2(b)} & \uparrow (n-1)1n\q, \mbox{ where }\q=4,5,\ldots,n-2         & (n-5)x^2                                \\  \hline
			\mbox{2(c)} & \uparrow (n-2)1n(n-1)                                       & x^2                                     \\ \hline
			\mbox{3(a)} & $\st n2\q$,\  \mbox{ where } \q\neq 1                       & (n-3)x + (n-3)(n-4) x^2                 \\
			\hline
			%	\mbox{3(b)} & \uparrow n12 & x \\
			%\hline
			\mbox{4(c)}
			            & \uparrow n1\q\q, \mbox{ where } \q\q \mbox{ is increasing }
			            & {n-2\choose 2}x                                                                                       \\[3pt] \hline
		\end{array}$\bigskip
	\caption{Permutations in $\exactly{n}{n-3}^{321}$ for {$n\geq5$}}\label{thelist-321}
\end{table}

\subsection{Proof of Theorem~\ref{213-thm}}

\begin{table}[h!]
	{		\renewcommand{\arraystretch}{1.05}
		$\begin{array}{|l|c|l|} \hline\hline
				                                        &                                                          &                                               \\[-2ex]
				\mbox{Case}                             & \mbox{Type of } w                                        & \mbox{Descent gen. polynomial }  A_n \glob{w} \\[-2ex]
				                                        &                                                          &                                               \\\hline\hline
				\mbox{2(a)}                             & \uparrow n3                                              & x                                             \\  \hline
				\mbox{3(a)}                             & \uparrow n23                                             & x                                             \\  \hline
				\mbox{3(b)}                             & 1\uparrow n\q2                                           & (n-3)x^2\quad (n-3\mbox{ choices for }\q)     \\  \hline
				\mbox{4(a)}                             & \uparrow n\q\q 1, \mbox{ but not } \uparrow n(n-2)(n-1)1 & (n-4)x^2+{n-2\choose 2}x^3                    \\  \hline
				\mbox{4(b)}                             & \uparrow n\q12                                           & (n-3)x^2 \quad (n-3\mbox{ choices for }\q)    \\  \hline
				\mbox{4(c)}                             & \uparrow n123  \mbox{ and } \uparrow n132                & x+x^2                                         \\  \hline
				\mbox{4(d)}                             & \uparrow n(n-2)(n-1)2                                    & x^2                                           \\  \hline
				\mbox{5(a)}                             & A n B (n-2)(n-1)2,\ AB=\uparrow,\ B \neq \emptyset       & (n-4)x^2                                      \\ \hline
				\mbox{5(b)}                             & A n B (n-3)(n-1)(n-2)1,\ AB=\uparrow                     & (n-4)x^3                                      \\ \hline
				\mbox{5(c)}                             & A n B (n-1) C (n-3)(n-2)1,\ ABC=\uparrow                 & {n-3\choose 2}x^3                             \\ \hline
				\mbox{5(d)}                             & A(n-1)nB(n-4)(n-3)(n-2)1                                 & (n-5)x^2                                      \\
				                                        & AB = \uparrow                                            &                                               \\ \hline
				\mbox{5(e)}                             & AnB (n-2) (n-1)\!\left\{\!\!\!
				\begin{array}{c}
					12 \\ \q 1
				\end{array}\!\!\!\right\},\ AB=\uparrow & (n-4)\left(x^2+(n-4)x^3\right)                                                                           \\  \hline
			\end{array}$\bigskip}
	\caption{Permutations in $\exactly{n}{n-3}^{213}$ for {$n\geq6$}}\label{thelist-213}
\end{table}

Considering (\ref{West-n-2-exact}) we can see that for $n\geq 4$,
\begin{align}\label{West-n-2-exact-213}
	\exactly{n}{n-2}^{213} = \sym_n \cap
	\glob{\;\uparrow n2,\; \uparrow n12,\; \uparrow n\q1,\; A n B (n-2)(n-1)1}
\end{align}
where to avoid the pattern 213, $AB=\uparrow$. Note that the pattern $\st(n-1)1n$ in $\exactly{n}{n-2}$ disappears because it contains the occurrence $(n-1)1n$ of the pattern 213.

It is easy to see that $A_n \glob{\uparrow n2}=x$, $A_n \glob{\uparrow n12}=x$, $A_n \glob{\uparrow n\q1}=(n-2)x^2$ (choosing $\q$ in $(n-2)$ ways) and $A_n \glob{\uparrow n\uparrow(n-2)(n-1)1}=(n-3)x^2$ (choosing placement of $n$ in $(n-3)$ ways), and hence from (\ref{West-n-2-exact-213}) we have
$$E^{213}_{n,n-2}(x)=2x+(2n-5)x^2.$$
Now, since $\exactly{n}{n-1}=\sym_n \cap \glob{\st n1}$ (\cite{West1990Permutations}),  $\exactly{n}{n-1}^{213}=\sym_n \cap \glob{\uparrow n1}$ and $E^{213}_{n,n-1}(x)=x$. Since $\sortable{n}{n-1}=\sym_n$, formula (\ref{1-st-sort-132}) can be used to find $W^{213}_{n,n-1}(x)$, and then iteration of (\ref{eq-W-E}) for $t=n-1, n-2$ gives the desired formulas for $W^{213}_{n,n-2}(x)$ and $W^{213}_{n,n-3}(x)$.

To complete the proof of Theorem~\ref{213-thm}, we need to compute $E^{213}_{n,n-3}(x)$ and apply (\ref{eq-W-E}) for $t=n-3$. We analyse Tables~\ref{thelist1} and~\ref{thelist2} presenting $\exactly{n}{n-3}$ using the observation that to the left of a large element, smaller elements must be in increasing order to avoid the pattern 213. We see that to avoid the pattern 213, the following cases are not possible: 1(a)--1(e), 2(b), 2(c), 5(f)--5(h) (as well as the first subcase of 5(d)). The remaining cases are listed in Table~\ref{thelist-213} and they give the desired formula for $W^{213}_{n,n-4}(x)$. Here are our comments for some cases in Table~\ref{thelist-213}.
\begin{itemize}
	\item[4(a)]  Suppose $a<b$. Then two possibilities are $\uparrow nab1$ and $\uparrow nba1$. Note that in the latter case any choice of $a,b$ gives a 213-avoiding permutation (justifying the term of ${n-2\choose 2}x^3$), while in the former case to avoid $213$, we must have $b=a+1$. Since in the former case choosing $a=n-2$ is not an option by the 4(a) description, $a$ can be chosen from $\{2,3,\ldots, n-3\}$, justifying the term of $(n-4)x^2$.
	\item[5(a)]  $(n-4)$ is the number of ways to insert $n$ in $AB$ of length $n-4$ so that $B$ is non-empty.
	\item[5(c)]  ${n-3\choose 2}$ is the number of ways to insert $n$ and $n-1$ in $ABC$ of length $n-5$. Any insertion will result in 3 descents in the resulting permutation.
	\item[5(d)]  $(n-5)$ is the number of ways to insert $(n-1)n$ in $AB$ of length $n-6$. Any insertion  will  result in 2 descents in the outcome permutation.
	\item[5(e)]  In $AnB(n-2)(n-1)12$ we have 2 descents and $n$ can be placed in $(n-4)$ ways in the permutation $AB$ of length $n-5$. On the other hand, independently to the choices of placing $n$ in  $AnB(n-2)(n-1)\q1$ with 3 descents, we can choose $\q$ in $(n-4)$ ways ({any choice will result in a 213-avoiding permutation}).
\end{itemize}

\subsection{Proof of Theorem~\ref{132-thm}}

Considering (\ref{West-n-2-exact}) we can see that for $n\geq 4$,
\begin{align}\label{West-n-2-exact-132}
	\exactly{n}{n-2}^{132} = \sym_n \cap
	\glob{\; \st(n-1)1n,\; \st n12,\; \st n21,\; n\st(n-2)\st (n-1)1\;}
\end{align}
where $\st$ represents any $132$-avoiding permutation in the first three patterns, and {$\st(n-2)\st$ represents any non-empty $132$-avoiding permutation} on elements $\{2,3,\ldots,n-2\}$. Indeed, for example, in the unrestricted set $\exactly{n}{n-2}$, the pattern $\st n2$ is not taken to the set  $\exactly{n}{n-2}^{132}$ as any permutation corresponding to it contains the pattern 132 (formed by the elements in $\{1,2,n\}$). Also, the leftmost $\st$ in $\st n\st(n-2)\st (n-1)1$ must be empty as otherwise any element in it along with $n$ and $n-1$ will form an occurrence of the pattern 132. Moreover, $\q=2$ in $\st n1\q$ or otherwise the element in place of $\q$ along with $2$ and $n$ will form an occurrence of the pattern 132.

It is easy to see that $A_n \glob{\st(n-1)1n}=xN_{n-3}(x)$, $A_n \glob{\st n12}=xN_{n-3}(x)$, $A_n \glob{\st n21}=x^2N_{n-3}(x)$ and, as $n\geq 4$, $A_n \glob{n\st(n-2)\st (n-1)1}=x^2N_{n-3}(x)$. Hence, it follows from \eqref{West-n-2-exact-132} that
$$E^{132}_{n,n-2}(x)=2(x+x^2)N_{n-3}(x).$$
Now, since $\exactly{n}{n-1}=\sym_n \cap \glob{\st n1}$ (\cite{West1990Permutations}),  $\exactly{n}{n-1}^{132}=\glob{\st n1}$, where $\st$ is any 132-avoiding permutation on $\{2,3,\ldots,n-1\}$, and $E^{132}_{n,n-1}(x)=xN_{n-2}(x)$. Since $\sortable{n}{n-1}=\sym_n$, formula (\ref{1-st-sort-132}) can be used to find $W^{132}_{n,n-1}(x)$, and then iteration of (\ref{eq-W-E}) for $t=n-1, n-2$ gives the desired formulas for $W^{132}_{n,n-2}(x)$ and $W^{132}_{n,n-3}(x)$.

\begin{table}[h!]
	{		\renewcommand{\arraystretch}{1.05}
		$\begin{array}{|l|c|l|} \hline\hline
				                          &                                      &                                               \\[-2ex]
				\mbox{Case}               & \mbox{Type of } w                    & \mbox{Descent gen. polynomial }  A_n \glob{w} \\[-2ex]
				                          &                                      &                                               \\\hline\hline
				\mbox{1(b)}               & \st(n-1)21n                          & x^2N_{n-4}(x)                                 \\   \hline
				\mbox{1(c)}               & \st(n-1)12n                          & xN_{n-4}(x)                                   \\   \hline
				\mbox{1(d)}               & (n-1)\st (n-3) \st (n-2)1n           & x^2N_{n-4}(x)                                 \\   \hline
				\mbox{1(e)}               & \st(n-2)1(n-1)n                      & xN_{n-4}(x)                                   \\   \hline
				%			\mbox{2(c)} & \st(n-2)1n(n-1) & x^2N_{n-4}(x) \\   \hline
				%			\mbox{3(a)} & \st n21&   x^2N_{n-3}(x) \\ \hline
				%			\mbox{3(b)} & \st n12 &  xN_{n-3}(x) \\ \hline
				\mbox{4(a)}               & \st n231 \mbox { and } \st n321      &
				%			\begin{array}{c} x^3 \mbox{ if } n=4; \mbox{ otherwise } \\
				(x^2+x^3)N_{n-4}(x)                                                                                              \\   \hline
				\mbox{4(b)}               & \st n213 \mbox{ and } \st n312       & 2x^2 N_{n-4}(x)                               \\   \hline
				\mbox{4(c)}               & \st n123                             & x N_{n-4}(x)                                  \\   \hline
				%			\mbox{5(b)} & n\st (n-3)\st (n-1)(n-2)1 & x^3N_{n-4}(x)\\   \hline
				\mbox{5(c)}               & n(n-1)\st (n-3)\st(n-2)1             & x^3N_{n-4}(x)                                 \\   \hline
				\mbox{5(d)}               &
				(n-1)  n \st\!\left\{\!\!\!
				\begin{array}{c}
					(n-3) \st (n-4) \\
					(n-4) \st (n-3)
				\end{array}
				\!\!\!\right\}\!\st (n-2)1
				                          & x^2N_{n-4}(x)
				\\  \hline
				\mbox{5(e)}               & n \st (n-2) \st (n-1)\!\left\{\!\!\!
				\begin{array}{c}
					12 \\ 2 1
				\end{array}\!\!\!\right\} & (x^2+x^3)N_{n-4}(x)                                                                  \\   \hline
				%			\mbox{5(f)} & 52413 & x^2 \mbox{ for } n=5 \\   \hline
				\mbox{5(g)}               & n\st (n-3)\st (n-2)1(n-1)            & x^2N_{n-4}(x)                                 \\   \hline
			\end{array}$\bigskip}
	\caption{Permutations in $\exactly{n}{n-3}^{132}$ for $n\geq 6$}\label{thelist-132}
\end{table}

To complete the proof of Theorem~\ref{132-thm}, we need to compute $E^{132}_{n,n-3}(x)$ and apply (\ref{eq-W-E}) for $t=n-3$. We analyse Tables~\ref{thelist1} and~\ref{thelist2} presenting $\exactly{n}{n-3}$. We see that to avoid the pattern 132, none of the following cases are  possible: 1(a), 2(a), 2(b), 2(c), 3(a), 3(b), 4(d), 5(a), 5(b), 5(f), 5(h); in particular, in 4(d) the elements in $\{1,2,n\}$ form the pattern 132. The remaining cases are listed in Table~\ref{thelist-132} and they give the desired formula for $W^{132}_{n,n-4}(x)$. We comment three cases in Table~\ref{thelist-132}.
\begin{itemize}
	\item[1(b)] $x^2$ records the descents $(n-1)2$ and 21, and {$\st$ is any $132$-avoiding permutation} on $\{3,4,\ldots,n-2\}$.
	\item[1(d)]  No matter if $\st$ is empty or not in $(n-1)\st(n-3)$, $(n-1)$ will be involved in a descent, and $(n-2)1$ is the other descent recorded by $x^2$. Also, $\st(n-3)\st$ is any non-empty 132-avoiding permutation on $\{2,3,\ldots,n-3\}$.
	\item[5(d)]   Note that $\st(n-3)\st(n-4)\st$ or $\st(n-4)\st(n-3)\st$ can be any non-empty 132-avoiding permutation on $\{2,3,\ldots,n-3\}$. This gives the term of $x^2N_{n-4}(x)$.
	      %Moreover, it is easy to see that
	      %$$A_n \glob{(n-1)n\st(n-4) \st (n-3)\st(n-2)1}=x^2A_{n-4} \glob{\st(n-4) \st (n-3)\st}=$$
	      %$$x^2\left(A_{n-4} \glob{\st(n-4)(n-3)\st}+A_{n-4} \glob{\st(n-4) \q (n-3)\st}\right)=$$
	      %$$x^2\left(W^{132}_{n-5,n-6}(x)+A_{n-4} \glob{\st(n-4)1(n-3)}+A_{n-4} \glob{\st(n-4) \q (n-3)\st1\st}\right)$$
	      %$$=x^2\left(W^{132}_{n-5,n-6}(x)+xW^{132}_{n-7,n-8}(x)+x^2\sum_{i=0}^{n-4}W^{132}_{i,i-1}(x)W^{132}_{n-i-7,n-i-8}(x)\right)$$
	      %where in $\st(n-4) \q (n-3)\st1\st$, $\st1\st$ is any 132-avoiding permutation formed by the smallest elements and $\q$ is the next smallest element, while the leftmost $\st$ is any 132-avoiding permutation formed by the largest elements; the sum is then obtained by letting the size of the leftmost $\st$ be $i$.
\end{itemize}

\subsection{Proof of Theorem~\ref{312-thm}}

Considering (\ref{West-n-2-exact}) we can see that for $n\geq 4$,
\begin{align}\label{West-n-2-exact-312}
	\exactly{n}{n-2}^{312} = \sym_n \cap
	\glob{\;1\st n2,\; \st (n-1)1n,\; \st n\q1}
	%         ,\; \st n(n-2)(n-1)1}
\end{align}
where $\st$ can be any $312$-avoiding permutation over the respective set.  Note that the pattern $\st n1\q$ in $\exactly{n}{n-2}$ disappears because it contains the occurrence $n1\q$ of the pattern 312.

It is easy to see that $A_n \glob{1\st n2}=A_n \glob{\st (n-1)1n}=xN_{n-3}(x)$. Computing $A_n\glob{\st n\q1}$ requires the following lemma.

\begin{lem}\label{312-end-asc} {For $n\geq2$,} the descent generating polynomial for $312$-avoiding permutations of length $n$ that end with
	\begin{itemize}
		\item  an ascent is $N_{n-1}(x)$;
		\item a descent is $N_n(x)-N_{n-1}(x)$.
	\end{itemize}
\end{lem}

\begin{proof} To avoid the pattern $312$, any permutation ending with an ascent must end with the largest element $n$. But then the element $n$ does not contribute any descent and is independent from the rest of the permutation hence proving the first claim. The second claim is straightforward as any permutation either ends with an ascent or with a descent.\end{proof}

Now, $A_n\glob{\st n\q1}=xA_{n-1}\glob{\st (n-1)\q}$ where $\st\q$ can be any $312$-avoiding permutation of length $(n-2)$. If $\st\q$ ends with an ascent, an extra descent will be created after inserting $(n-1)$, so by Lemma~\ref{312-end-asc} descents are counted by $xN_{n-3}(x)$. If $\st\q$ ends with a descent, inserting $(n-1)$ does not create an extra descent, so by Lemma~\ref{312-end-asc} descents are counted by $N_{n-2}(x)-N_{n-3}(x)$. Hence
\begin{equation}\label{stnq1}
	A_n\glob{\st n\q1}=xN_{n-2}(x)+x(x-1)N_{n-3}(x),
\end{equation} and from (\ref{West-n-2-exact-312}) we have
$$E^{312}_{n,n-2}(x)=xN_{n-2}(x)+x(x+1)N_{n-3}(x).$$
Now, since $\exactly{n}{n-1}=\sym_n \cap \glob{\st n1}$ (\cite{West1990Permutations}),  $\exactly{n}{n-1}^{312}=\sym_n \cap \glob{\st n1}$ and $E^{312}_{n,n-1}(x)=xN_{n-2}(x)$. Since $\sortable{n}{n-1}=\sym_n$, formula (\ref{1-st-sort-132}) can be used to find $W^{312}_{n,n-1}(x)$, and then iteration of (\ref{eq-W-E}) for $t=n-1, n-2$ gives the desired formulas for $W^{312}_{n,n-2}(x)$ and $W^{312}_{n,n-3}(x)$.

To complete the proof of Theorem~\ref{312-thm}, we need to compute $E^{312}_{n,n-3}(x)$ and apply (\ref{eq-W-E}) for $t=n-3$. We analyse Tables~\ref{thelist1} and~\ref{thelist2} presenting $\exactly{n}{n-3}$ using the observation that to the right of a large element, smaller elements must be in decreasing order to avoid the pattern 312. We see that to avoid the pattern 312, the following cases are not possible: 1(c),1(d), 2(b), 3(a), 4(b)--4(d), 5(a-h), and also the restriction ``but not ...'' in 4(a) can be removed because $n(n-2)(n-1)$ forms the pattern 312. The remaining cases are listed in Table~\ref{thelist-312} ($\st$ there denotes any $312$-avoiding permutation over the respective set) and they give the desired formula for $W^{312}_{n,n-4}(x)$. We comment on three cases in Table~\ref{thelist-312}.
\begin{itemize}
	\item[1(b)]   {$A_n\glob{\st (n-1)\q 1n}=A_{n-1}\glob{\st (n-1)\q1}$} and is given by (\ref{stnq1}).
	\item[3(b)]   $A_n\glob{1\st n\q 2}=A_{n-1}\glob{\st (n-1)\q1}$ and is given by (\ref{stnq1}).
	\item[4(a)]  {Here the notation denotes the union over all $b>a$. After deleting the final $1$, standardizing, and relabeling the two variables,}
	\[
	{A_n\glob{\st nba1}=xA_{n-1}\glob{\st (n-1)ba}.}
	\]
	To compute $A_{n-1}\glob{\st (n-1)ba}$, note that the elements $1,2,\ldots,a-1$ in $\st$ must be in the leftmost $a-1$ positions (or else there will be an occurrence of the pattern 312). Hence, the permutation formed by the smallest $a-1$ elements is independent from the rest of a 312-avoiding permutation $\pi$ of length $n-1$, and
	\begingroup
	\begin{align*}
	A_{n-1}\glob{\st (n-1)ba}
	=\sum_{a=1}^{n-3}N_{a-1}(x)
	\left(x(x-1)N_{n-a-3}(x)+xN_{n-a-2}(x)\right).
	\end{align*}
	\endgroup
	Here $\st(n-a)\q1$ is the pattern formed by the $n-a$ largest elements in $\pi$, and formula~(\ref{stnq1}) is applied.

	      By using the following recurrence relation of Narayana polynomials~(see~\cite[Theorem 2.2]{Petersen2015})
	      \begin{align*}
		      {N_n(x)=N_{n-1}(x)+x\sum_{i=0}^{n-2}N_i(x)N_{n-1-i}(x)} \mbox{   for }  n \geq 1,
	      \end{align*}
	      we get that
	      \begin{align*}
		      A_{n-1}\glob{\st (n-1)ba}
		       & = {(x-1)\left( N_{n-3}(x) +(x-1)N_{n-4}(x)  \right)} +  \left( N_{n-2}(x)-N_{n-3}(x) \right) \\
		       & = N_{n-2}(x) + (x-2) N_{n-3}(x) +  (x-1)^2 N_{n-4}(x).
	      \end{align*}

\end{itemize}

\begin{table}[h!]
	%{		\renewcommand{\arraystretch}{1.05}
	$\begin{array}{|l|c|l|} \hline\hline
			            &                                &                                                            \\[-2ex]
			\mbox{Case} & \mbox{Type of } w              & \mbox{Descent gen. polynomial }  A_n \glob{w}              \\[-2ex]
			            &                                &                                                            \\\hline\hline
			\mbox{1(a)} & 1\st (n-1)2n                   & xN_{n-4}(x)                                                \\  \hline
			\mbox{1(b)} & \st(n-1)\q1n                   & xN_{n-3}(x) + x(x-1)N_{n-4}(x)                             \\  \hline
			\mbox{1(e)} & \st (n-2)1(n-1)n               & xN_{n-4}(x)                                                \\  \hline
			\mbox{2(a)} & 12\st n3 \mbox{ and } 21\st n3 & x(x+1)N_{n-4}(x)                                           \\  \hline
			\mbox{2(c)} & \st(n-2)1n(n-1)                & x^2N_{n-4}(x)                                              \\  \hline
			%	\mbox{3(a)} & \st n21 & x^2N_{n-3}(x)  \\  \hline
			\mbox{3(b)} & 1\st n\q 2                     & xN_{n-3}(x) + x(x-1)N_{n-4}(x)                             \\  \hline
			\mbox{4(a)} & \st nba1,\ b>a                 & x  N_{n-2}(x) + x (x-2) N_{n-3}(x) +  x (x-1)^2 N_{n-4}(x) \\  \hline
			%	\mbox{5(c)} &  \st n(n-1)(n-3)(n-2)1 & x^3W^{312}_{n-5,n-6}(x),\ n\geq 5\\  \hline
		\end{array}$\bigskip
	\caption{Permutations in $\exactly{n}{n-3}^{312}$ for {$n\geq5$}}\label{thelist-312}
\end{table}

\section{Concluding remarks}\label{final-sect}

Table~\ref{length-3-pattern-avoidance} confirms the following conjecture in the cases of $t\in\{n-1,n-2,n-3,n-4\}$, while the case of $t=1$ is equivalent to avoiding two patterns (additionally the pattern 231 in each case) and is known, and is easy to see to be true \cite{Kitaev2011Patterns}.

\begin{conj}\label{213-321-conj} The number of $213$-avoiding $t$-stack-sortable permutations of length $n$ is the same as that of $321$-avoiding  $t$-stack-sortable permutations of length $n$ for any $t\geq 0$ and $n\geq 0$.  \end{conj}

Conjecture~\ref{213-321-conj} has also been confirmed computationally for  $t\in\{2,3,4,5\}$ for small values of $n$, and the respective numbers, along with those corresponding to $t=1$, are given in Table~\ref{table-comput-res}. Interestingly, Table~\ref{table-comput-res} matches the table in Section~4 in \cite{Cerbai2020Stack}, which leads to the following conjecture related to so-called 321-machine that we state in terms of pattern avoidance (again, based on  Section~4 in \cite{Cerbai2020Stack}).

\begin{table}
	\begin{center}
		\begin{tabular}{|c|l|}
			\hline
			$t=1$ & 1, 2, 4, 8, 16, 32, 64, 128, 256, 512, 1024, $\ldots$         \\
			\hline
			$t=2$ & 1, 2, 5, 13, 34, 89, 233, 610, 1597, 4181, 10946, $\ldots$    \\
			\hline
			$t=3$ & 1, 2, 5, 14, 41, 122, 365, 1094, 3281, 9842, 29525, $\ldots$  \\
			\hline
			$t=4$ & 1, 2, 5, 14, 42, 131, 417, 1341, 4334, 14041, 45542, $\ldots$ \\
			\hline
			$t=5$ & 1, 2, 5, 14, 42, 132, 428, 1416, 4744, 16016, 54320, $\ldots$ \\
			\hline
		\end{tabular}
		\caption{The number of $t$-stack-sortable $p$-avoiding permutations of length $n\geq 1$ for $p\in\{213, 321\}$}\label{table-comput-res}
	\end{center}
\end{table}

\begin{conj}\label{conj2} For any $t\geq 1$, $t$-stack-sortable $p$-avoiding permutations of length $n\geq 0$, for $p\in\{213, 321\}$, are in one-to-one correspondence with permutations of length $n$ avoiding the patterns $132$ and $12\cdots (t+2)$ simultaneously. \end{conj}

Research directions for $t$-stack-sortable permutations include the  study of unimodality, log-concavity and real-rootedness  of the polynomials $W_{n,t}(x)$ for any $1 \le t \le n-1$.
We will not define these notions here, instead referring the interested readers to ~\cite{Bona2002Symmetry,Braenden2008Actions,Braenden2015Unimodality,Brenti1994Log,Stanley1989Log}. In particular, B\'ona~\cite{Bona2002Symmetry} proved  the symmetry and unimodality of $W_{n,t}(x)$ by giving combinatoiral bijecitons and injections.  Based on numerical evidence, B\'ona further conjectured log-concavity and real-rootedness of $W_{n,t}(x)$.
Br{\"a}nd{\'e}n~\cite{Braenden2008Actions} strengthened B\'ona's result by showing its $\gamma$-positivity via a combinatorial group action.
Since no exact formulas are known for  $W_{n,t}(x)$ in general, the remaining log-concavity and real-rootedness problems seem to be very hard. The real-rootedness of Narayana polynomials and Eulerian polynomials are well-known~\cite{Liu2007unified} and the real-rootedness of $W_{n,n-2}(x)$ was  proved by  Br{\"a}nd{\'e}n~\cite{Braenden2006linear} and Zhang~\cite{Zhang2015real}.

In this paper,  we obtained expressions for $W_{n,t}(x)$ in the cases of $t=n-3$ and $t=n-4$ in terms of Eulerian polynomials. There are several known approaches to study the unimodality, $\gamma$-positivity,  log-concavity and  real-rootedness for Eulerian polynomials, see \cite{Braenden2008Actions, Braenden2015Unimodality} for example. In particular, Br{\"a}nd{\'e}n~\cite{Braenden2008Actions} proved combinatorially $\gamma$-positivity for $W_{n,t}(x)$ for any $t$ that also implies unimodality of $W_{n,t}(x)$. However, we still cannot prove the  log-concavity and  real-rootedness for $W_{n,t}(x)$. The formulas derived in this paper could potentially be useful to answer the questions, for example, if more analytic properties of Eulerian polynomials and/or the stack-sorting operation are discovered.

Finally, $p$-avoiding $t$-stack-sortable permutations can be studied for longer patterns $p$, or indeed for other types of patterns $p$ ({\em consecutive}, {\em vincular}, {\em bivincular}, etc \cite{Kitaev2011Patterns}, not to be defined here). Such studies are interesting in their own right, but also may bring interesting connections to other combinatorial objects.

\vskip 3mm
\section*{Acknowledgments}
Philip B. Zhang was supported by the National Natural Science Foundation of China (No. 12171362).
The authors thank Larry X. W. Wang for bringing the incorrect formula in Theorem~\ref{main-thm2} to our attention.

%\section*{Declarations}
%\textbf{Conflict of Interest}  On behalf of all authors, the corresponding author states that there is no conflict of interest.

\end{document}